\documentclass[11pt]{article}

\usepackage{a4wide}
\usepackage{amsmath,amsbsy,amssymb,latexsym,amsthm}

\newtheorem{Teo}{Theorem}
\newtheorem{Prop}{Proposition}
\newtheorem{Remark}{Remark}

\renewenvironment{abstract}{\begin{center}
\begin{minipage}[c]{12cm}
{\begin{center}\bf Abstract \end{center}}}
{\end{minipage}
\end{center}}
\newenvironment{keywords}{\begin{center}
\begin{minipage}[c]{12cm}
{\bf Keywords:}} {\end{minipage}
\end{center}}
\newenvironment{msc}{\begin{center}
\begin{minipage}[c]{12cm}
{\bf Mathematics Subject Classification 2010:}}
{\end{minipage}
\end{center}}

\begin{document}

\title{Leitmann's direct method
for fractional optimization problems\thanks{Dedicated to
George Leitmann on the occasion of his 85th birthday.
The authors are grateful to George Leitmann
for the suggestions regarding improvement of the text.
Submitted June 16, 2009 and accepted March 15, 2010
for publication in \emph{Applied Mathematics and Computation}.}}

\author{Ricardo Almeida\\
        \texttt{ricardo.almeida@ua.pt}
   \and Delfim F. M. Torres\\
        \texttt{delfim@ua.pt}}

\date{Department of Mathematics\\
University of Aveiro\\
3810-193 Aveiro, Portugal}

\maketitle


\begin{abstract}
Based on a method introduced by Leitmann
[Internat. J. Non-Linear Mech. {\bf 2} (1967), 55--59],
we exhibit exact solutions for some fractional optimization
problems of the calculus of variations and optimal control.
\end{abstract}

\begin{msc}
26A33, 49J05, 49J15.
\end{msc}

\begin{keywords}
calculus of variations, optimal control, direct methods,
fractional calculus, Leitmann's principle.
\end{keywords}


\section{Introduction}
\label{sec1}

A direct method to determine solutions for some problems
of the calculus of variations, avoiding the variational approach,
was introduced by George Leitmann in 1967 \cite{Leitmann}.
Leitmann's method is based on the following two steps:
(i) by an appropriate coordinate transformation,
we rewrite the initial problem as an equivalent simpler one;
(ii) knowing the solution for the new equivalent problem,
and since there exists an one-to-one correspondence between
the minimizers (or maximizers) of the new problem with
the ones of the original problem, we determine the desired solution.
We also note \cite{silva,torres}, where this method is extended
to optimal control. The idea behind the direct method
of Leitmann is clever and simple,
providing a ``principle'' --- a general assertion
that can be formulated in precise terms, as a theorem,
in many different contexts, under many different assumptions.
Leitmann's principle reflects the fact that optimality
does not depend on the particular coordinates used,
being preserved under invariance transformations
up to an exact differential.
For more on the subject we refer the reader to
\cite{Carlson,Leitmann2,Leitmann3,torres}.

Fractional calculus is a branch of mathematics where one allows integrals
and derivatives to be of an arbitrary real or complex order \cite{miller,Podlubny,book:Samko}.
It has found many applications in pure and applied mathematics,
physics, and engineering \cite{Debnath,kilbas}. An area of fractional calculus
that is receiving remarkable attention consists in developing a theory of the calculus
of variations and optimal control in the presence
of fractional (non-integer) derivatives, and investigating its applications
--- see \cite{AGRA1,Baleanu:D:A,Cresson,El-Nabulsi:Torres:2008,Frederico:Torres:2008}
and references therein.

The literature on necessary optimality conditions for
fractional optimal control problems is vast (\textrm{cf.}, \textrm{e.g.},
\cite{AGRA1,ric:delf,Baleanu,Baleanu:D:A,El-Nabulsi:Torres:2007,Frederico:Torres:2007}).
Applications of such necessary conditions involve, however, the solution of
fractional differential equations, which is a difficult task
to accomplish analytically. Often, numerical methods or truncated Taylor series
are used \cite{AGRA2,Baleanu:D:A,Diethelm}. In this note,
motivated by the results of \cite{withGasta:SI:Leit:85,silva,torres},
we give the first examples of fractional problems of the calculus
of variations and optimal control with a known
exact expression for the global minimizers.
Such minimizers are found utilizing Leitmann's principle.

Although simple and easy to understand,
Leitmann's principle is not sufficiently well-known and used.
This is particularly true in the present context
of fractional calculus, where Leitmann's
direct method went completely unnoticed while
necessary optimality conditions, extremely difficult
to apply and often proved
under (too) strong hypotheses, are widely explored.
Recently there has been an increase of interest
in Leitmann's method with several papers devoted
to various aspects of the subject,
including discrete-time problems
of the calculus of variations and, more generally,
variational problems on time scales \cite{abmalina:delfim}.
There is, however, no direct link between
the time scale calculus and the fractional calculus:
they are different areas of research, with
different communities. Moreover,
there are profound and interesting implications of Leitmann's
method for fractional problems (\textrm{cf. e.g.} Remark~\ref{rem:resp:r2})
that have no counterpart in other problems of the calculus
of variations, such as the ones studied in
\cite{Ferreira:Torres:2007,abmalina:delfim,natalia:delfim}.
This is due to the fact that Euler-Lagrange
type conditions proved in the literature
of the fractional calculus of variations are
often deduced imposing that admissible functions
have continuous left and right fractional derivatives
on the closed interval where the problem is defined
--- see \cite{AGRA1} and discussions in \cite{MyID:182,ric:delf}.
This is different from the calculus of variations in different contexts,
\textrm{e.g.}, different from what happens with time scale variational problems
\cite{abmalina:delfim}.

The text is organized as follows. In Section~\ref{sec2} we present a short introduction
to fractional calculus, providing the necessary concepts and results needed in the sequel.
A simple version of Leitmann's method for solving fractional variational problems
is given in Section~\ref{sec3}, while in Section~\ref{sec4}
we illustrate the effectiveness of the method
by obtaining exact global minimizers for fractional variational problems depending
on fractional derivatives and integrals,
as well as fractional optimal control problems.
We end with conclusions in Section~\ref{sec:conc}.


\section{Preliminaries}
\label{sec2}

Let $f:[a,b]\rightarrow\mathbb{R}$ be a continuous function,
$\alpha > 0$ a given real, and $n:=[\alpha]+1$, where $[\alpha]$
denotes the integer part of $\alpha$. Let $\Gamma(\cdot)$
be the Gamma function, \textrm{i.e.},
$$
\Gamma(z)=\int_0^\infty t^{z-1}e^{-t}\, dt, \quad z>0.
$$
The Riemann-Liouville fractional integral of order $\alpha$ is defined by
$$
{_aI_x^\alpha}f(x)=\frac{1}{\Gamma(\alpha)}\int_a^x (x-t)^{\alpha-1}f(t)dt \, ,
$$
while the Riemann-Liouville fractional derivative of order $\alpha$ is given by
$$
{_aD_x^\alpha}f(x)=\frac{d^n}{dx^n}  {_aI_x^{n-\alpha}}f(x)
=\frac{1}{\Gamma(n-\alpha)}\frac{d^n}{dx^n}\int_a^x
(x-t)^{n-\alpha-1}f(t)dt \, .
$$
In particular,
$$
{_aD_x^0}f(x)=f(x) \quad \mbox{and}
\quad  {_aD_x^m}f(x)=f^{(m)}(x), \quad m \in \mathbb N.
$$

There exist composition formulas between fractional integrals
and fractional derivatives, similar to the standard ones
(proofs may be found, \textrm{e.g.}, in \cite{kilbas}):

\begin{Teo} If $f \in L_p(a,b)\,(1\leq p \leq \infty)$, $\alpha>0$ and $\beta>0$, then
$${_aI_x^\alpha}{_aI_x^\beta}f(x)={_aI_x^{\alpha+\beta}}f(x)$$
almost everywhere on $[a,b]$.
\end{Teo}

\begin{Teo} If $f \in L_p(a,b)\,(1\leq p \leq \infty)$, then
$${_aD_x^\alpha}{_aI_x^\alpha}f(x)=f(x)$$
almost everywhere on $[a,b]$.
\end{Teo}

\begin{Teo} Let $f_{n-\alpha}(x)={_aI_x^{n-\alpha}}f(x)$.
If $f \in L_1(a,b)$ and $f_{n-\alpha} \in AC^n[a,b]$, then
$${_aI_x^\alpha}{_aD_x^\alpha}f(x)=f(x)
-\sum_{j=1}^{n}\frac{f^{(n-j)}_{n
-\alpha}(a)}{\Gamma(\alpha-j+1)}(x-a)^{\alpha-j}$$
almost everywhere on $[a,b]$.
\end{Teo}

In particular, if ${_aI_a^{1-\alpha}}f(a)=0$
and $\alpha \in (0,1)$, then
${_aI_x^\alpha}{_aD_x^\alpha}f(x)=f(x)$.


\section{Fractional version of Leitmann's approach}
\label{sec3}

Consider the following fractional problem
of the calculus of variations: to minimize
\begin{equation}
\label{funct}
\mathfrak{J}(y):=\int_a^b F(x,y,{_aD^\alpha_x}y(x))\,dx
\end{equation}
under the constraint
\begin{equation}
\label{inte}
{_aI_b^{1-\alpha}}y(b)=y_b,
\end{equation}
where $y_b$ is a fixed real. Here we consider
continuous functions $y:[a,b]\to\mathbb R$ such that
${_aD^\alpha_x}y(x)$ exists on $[a,b]$. We also
emphasize that the condition ${_aI_a^{1-\alpha}}y(a)=0$
appears implicitly, since function $y$ is continuous
(\textrm{cf.} \cite[pag.~46]{miller}). If we allow
$\alpha=1$, functional (\ref{funct}) becomes the standard one,
\begin{equation}
\label{eq:class:funct}
\mathfrak{J}(y):=\int_a^b F(x,y,y'(x))\,dx \, ,
\end{equation}
and equality (\ref{inte}) reads as $y(b)=y_b$, that is,
we obtain a classical problem of the calculus of variations:
to minimize (\ref{eq:class:funct}) under given boundary
conditions $y(a) = 0$ and $y(b)=y_b$.

The next result is obtained by following
the method presented in \cite{Leitmann}.

\begin{Teo}
\label{thm:mr}
Let $y(x)=z(x,\tilde y (x))$ be a continuous
transformation having a unique inverse
$\tilde y(x)=\tilde z (x,y(x))$, such that
there exists an one-to-one correspondence
$$y(x) \Leftrightarrow \tilde y(x),$$
for every function $y$ satisfying
\begin{equation}
\label{constraint1}
{_aI_b^{1-\alpha}}y(b)=y_b
\end{equation}
and for every function $\tilde{y}$ satisfying
\begin{equation}
\label{constraint2}
{_aI_b^{1-\alpha}}\tilde y(b)={_aI_b^{1-\alpha}}\tilde z(x,y(x))|_{x=b}.
\end{equation}
In addition assume that there exists a $C^1$ function
$H:[a,b]\times\mathbb R\to\mathbb R$ such that the relation
$$
F(x,y,{_aD^\alpha_x}y(x))-F(x,\tilde y,{_aD^\alpha_x}\tilde y(x))
=\frac{d}{dx}H(x,{_aI^{1-\alpha}_x}\tilde y(x))
$$
holds. Then there exists a one-to-one correspondence
between the minimizers $y^{*}(x)$ of $\mathfrak{J}$
verifying (\ref{constraint1}) and the  minimizers
$\tilde{y}^*(x)=\tilde z (x,y^{*}(x))$ of $\mathfrak{J}$
verifying (\ref{constraint2}).
\end{Teo}

\begin{proof} This is obvious from
$$
\begin{array}{ll}
\mathfrak{J}({y^*})-\mathfrak{J}(\tilde{y}^*)
&=\displaystyle\int_a^b F(x,{y^*},{_aD^\alpha_x}{y^*}(x))
-  F(x,\tilde{y}^*,{_aD^\alpha_x}\tilde{y}^*(x))\,dx\\
&=\displaystyle\int_a^b \frac{d}{dx}H(x,{_aI^{1-\alpha}_x}\tilde y^*(x))\\
&=\displaystyle H(b,{_aI^{1-\alpha}_b}\tilde y^*(b))-H(a,0)
\end{array}$$
and the fact that the right-hand side is a constant.
\end{proof}

In the next section we show that Theorem~\ref{thm:mr} can be applied with success
not only to fractional variational problems of form (\ref{funct})--(\ref{inte})
but also to fractional optimal control problems \cite{AGRA2,Frederico:Torres:2008}
and problems depending on both fractional derivatives and fractional integrals
\cite{ric:delf}.


\section{Applications}
\label{sec4}

We illustrate the application of Leitmann's method
with examples from three different classes of problems:
Proposition~\ref{prop:prop} gives an exact solution
for a standard fractional variational problem (\ref{funct})--(\ref{inte});
Proposition~\ref{prop:to:ans:rv2} gives the global
minimizer for a family of problems whose Lagrangian
depends simultaneously on a fractional derivative
and a fractional integral of $y$;
Proposition~\ref{prop:2} gives the exact solution
of a fractional optimal control problem.

\begin{Prop}
\label{prop:prop}
Let $\alpha$ be an arbitrary real between 0 and 1.
The global minimizer of the fractional problem
\begin{equation}
\label{int1}
\mathfrak{J}(y)
=\int_0^1({_0D_x^\alpha}y(x))^2\, dx \longrightarrow \min
\end{equation}
under the constraint
\begin{equation}
\label{ex1:const}
{_0I_1^{1-\alpha}}y(1)=c
\end{equation}
is
\begin{equation}
\label{sol1}
y(x)=\frac{c}{\alpha\Gamma(\alpha)}x^\alpha.
\end{equation}
\end{Prop}

\begin{Remark}
We note that for $\alpha = 1$ problem (\ref{int1})--(\ref{ex1:const})
reduces to the standard problem of the calculus of variations
\begin{equation}
\label{ex1:class:form}
\begin{gathered}
\int_0^1 \left(y'(x)\right)^2\, dx \longrightarrow \min \\
y(0)=0 \, , \quad y(1)=c \, .
\end{gathered}
\end{equation}
The global minimizer of problem (\ref{ex1:class:form})
is $y(x)=cx$ (\textrm{cf. e.g.} \cite[pag.~509]{silva}),
which coincides with our solution (\ref{sol1}) for $\alpha=1$.
\end{Remark}

\begin{Remark}
\label{rem:resp:r2}
Problem (\ref{int1})--(\ref{ex1:const})
is trivial when $c = 0$:
since $\mathfrak{J}(y) \ge 0$ for any function $y$
and $\mathfrak{J}(0)$ is zero, it is obvious that the null function
is the solution of (\ref{int1})--(\ref{ex1:const}).
We are thus interested in the nontrivial case $c \ne 0$.
In \cite{AGRA1} a fractional Euler-Lagrange necessary
optimality condition is proved, under the hypothesis
that admissible functions $y$
have continuous left and right fractional derivatives
on the closed interval $[a,b]$. Such an assumption implies that
$y(a)=y(b)=0$ (\textrm{cf.} \cite{Ross:Samko:Love}).
This means that Agrawal's fractional Euler-Lagrange equation
given in \cite{AGRA1} can only be applied to (\ref{int1})--(\ref{ex1:const})
when we restrict ourselves to the trivial case $c=0$.
Indeed, if $c \ne 0$, then $y(1)=\frac{c}{\alpha\Gamma(\alpha)} \ne 0$
and therefore the right fractional derivative is not defined everywhere on
interval $[0,1]$. It is easy to check that for $c \ne 0$
the solution (\ref{sol1}) of problem (\ref{int1})--(\ref{ex1:const})
does not satisfy the Euler-Lagrange equation derived in \cite{AGRA1}.
The trivial situation $c = 0$ is the only case when
the methods of \cite{AGRA1} succeed in finding the solution
of (\ref{int1})--(\ref{ex1:const}).
\end{Remark}

\begin{proof}(of Proposition~\ref{prop:prop})
Let $y=\tilde{y}+f$, where $f$ is a function such that
${_0D_x^\alpha}f(x)=K$, for all $x \in [0,1]$
($K$ to be specified later). We consider functions
$\tilde{y}$ satisfying the equation
\begin{equation}
\label{system2}
{_0I_1^{1-\alpha}}\tilde y(x)=0.
\end{equation}
Since ${_0D_x^\alpha}f=K$, then
$$
{_0I_x^{\alpha}}\,{_0D_x^\alpha}f(x)={_0I_x^{\alpha}}K
$$
and the continuity on $f$ implies that ${_aI_a^{1-\alpha}}f(a)=0$.
Consequently
$$
f(x)=\frac{K}{\alpha\Gamma(\alpha)}x^\alpha.
$$
Moreover,
$$
\begin{array}{ll}
\displaystyle\int_0 ^1({_0D_x^\alpha}y)^2\, dx
& = \displaystyle\int_0^1({_0D_x^\alpha}\tilde y
+{_0D_x^\alpha}f)^2\, dx\\
&=\displaystyle \int_0^1({_0D_x^\alpha}\tilde y)^2 \, dx
+ \int_0^1(2K{_0D_x^\alpha}\tilde y+K^2)\, dx\\
&= \displaystyle\int_0^1({_0D_x^\alpha}\tilde y )^2\, dx
+ \int_0^1\frac{d}{dx}(2K{_0I_x^{1-\alpha}}\tilde y+K^2x)\, dx.
\end{array}
$$
Observe that $\int_0^1({_0D_x^\alpha}\tilde y)^2\geq0$
and $\int_0^1{_0D_x^\alpha}\tilde y \, dx=0$
if $\tilde{y}(x) \equiv 0$. Therefore $\tilde{y}(x) \equiv 0$
is a minimizer of $\mathfrak{J}$, and satisfies equation
(\ref{system2}). Since the second integral of the last sum
is constant (actually, is equal to $K^2$), it follows that
$$
y(x)=\frac{K}{\alpha\Gamma(\alpha)}x^\alpha
$$
is a solution of problem (\ref{int1})--(\ref{ex1:const}). Since
$$
{_0I_x^{1-\alpha}}\frac{K}{\alpha\Gamma(\alpha)}x^\alpha=Kx
$$
it follows that ${_0I_1^{1-\alpha}}y(1)=K$ and so $K=c$.
\end{proof}

We now study a variational problem
where not only a fractional derivative of
$y$ but also a fractional integral of $y$
appears in the Lagrangian. This
new class of fractional functionals of the calculus
of variations, that depend not only on fractional
derivatives but also on fractional integrals,
has been recently studied in \cite{ric:delf}.
In that paper necessary and sufficient conditions
of optimality for the fundamental problem of the
calculus of variations and for problems subject
to integral constraints (isoperimetric problems)
are studied. However, similarly to \cite{AGRA1},
the conditions are difficult to apply, and no example
has been solved analytically. We now show how Leitmann's
direct method can also be used to solve such problems.

\begin{Prop}
\label{prop:to:ans:rv2}
Let $g$ be a given function of class $C^1$ with $g(x)\ne 0$
on $[0,1]$, and $\alpha$ and $\xi$ real numbers with
$\alpha$ between zero and one.
The global minimizer of the fractional
variational problem
\begin{equation}
\label{(P)}
\begin{gathered}
\int_0^1\left[{_0D_x^\alpha}y(x)\cdot g(x)
+({_0I_x^{1-\alpha}}y(x)+1)g'(x)\right]^2\, dx
\longrightarrow \min\\
{_0I_1^{1-\alpha}}y(1)=\xi
\end{gathered}
\end{equation}
is given by the function
\begin{equation}
\label{eq:sol:prb_sem:paper}
y(x)={_0D_x^{1-\alpha}}\left( \frac{[g(1)(\xi+1)
-g(0)]x+g(0)}{g(x)}-1 \right).
\end{equation}
\end{Prop}

\begin{Remark}
For $\alpha = 1$ problem (\ref{(P)})
coincides with the classical problem
of the calculus variations
\begin{equation*}
\begin{gathered}
\int_0^1\left[y'(x)\cdot g(x)
+(y(x)+1)g'(x)\right]^2\, dx \longrightarrow \min\\
y(0) = 0\, , \quad y(1)=\xi \,
\end{gathered}
\end{equation*}
that has been studied by Leitmann
in his seminal paper of 1967 \cite{Leitmann}.
For the integer-order case $\alpha = 1$,
our function (\ref{eq:sol:prb_sem:paper})
reduces to $y(x)= \frac{[g(1)(\xi+1)-g(0)]x+g(0)}{g(x)}-1$,
which coincides with the results of Leitmann \cite{Leitmann}.
\end{Remark}

\begin{proof}(of Proposition~\ref{prop:to:ans:rv2})
To begin, observe that
$$
\frac{d}{dx}\left[ ({_0I_x^{1-\alpha}}y(x)+1)g(x) \right]
={_0D_x^\alpha}y(x)\cdot g(x)+({_0I_x^{1-\alpha}}y(x)+1)g'(x).
$$
Let $y$ be an admissible function to problem (\ref{(P)}),
and consider the transformation $y(x)=\tilde{y}(x) +f(x)$,
with $f$ to be determined later. Then
\begin{equation*}
\begin{split}
\Biggl(\frac{d}{dx} & \left[ ({_0I_x^{1-\alpha}}\tilde y(x)
+{_0I_x^{1-\alpha}}f(x)+1)g(x) \right]\Biggr)^2
-\left(\frac{d}{dx}\left[ ({_0I_x^{1-\alpha}}\tilde y(x)
+1)g(x) \right]\right)^2\\
&=2\frac{d}{dx}\left[ ({_0I_x^{1-\alpha}}\tilde y(x)+1)g(x) \right]
\cdot \frac{d}{dx}\left[ {_0I_x^{1-\alpha}}f(x)\cdot g(x) \right]
+\left(  \frac{d}{dx}\left[ {_0I_x^{1-\alpha}}f(x)
\cdot g(x) \right] \right)^2\\
&= \frac{d}{dx}\left[ {_0I_x^{1-\alpha}}f(x)\cdot g(x) \right]
\cdot \frac{d}{dx}\left[\left( 2({_0I_x^{1-\alpha}}\tilde y(x)+1)
+ {_0I_x^{1-\alpha}}f(x)\right)g(x) \right].
\end{split}
\end{equation*}
Let us determine $f$ in such a way that
$$
\frac{d}{dx}\left[ {_0I_x^{1-\alpha}}f(x)\cdot g(x) \right]=const.
$$
Integrating, we deduce that
$$
{_0I_x^{1-\alpha}}f(x)\cdot g(x)=Ax+B \, ,
$$
\textrm{i.e.},
$$
f(x)={_0D_x^{1-\alpha}}\left( \frac{Ax+B}{g(x)}\right).
$$
Consider now the new problem
\begin{equation}
\label{(P2)}
\begin{gathered}
\int_0^1\left({_0D_x^\alpha}\tilde y(x)\cdot g(x)
+({_0I_x^{1-\alpha}}\tilde y(x)+1)g'(x)\right)^2\,
dx \longrightarrow \min\\
{_0I_1^{1-\alpha}}\tilde y(1)=\frac{1}{g(1)}-1 .
\end{gathered}
\end{equation}
It is easy to see, and a trivial exercise to check, that
$$
\tilde{y}(x) = {_0D_x^{1-\alpha}}\left( \frac{1}{g(x)}-1 \right)
$$
is a solution of problem (\ref{(P2)}). Therefore
$$
y(x)={_0D_x^{1-\alpha}}\left( \frac{Ax+C}{g(x)}-1 \right),
\quad \mbox{with }C=B+1,
$$
is a solution of (\ref{(P)}).
Using the boundary conditions
$$
{_0I_0^{1-\alpha}}y(0)=0 \quad \mbox{and}
\quad{_0I_1^{1-\alpha}}y(1)=\xi
$$
we obtain the values for the constants $A$ and $C$:
$A=g(1)(\xi+1)-g(0)$ and $C=g(0)$.
\end{proof}

Proposition~\ref{prop:2} deals with a fractional optimal control problem.
Using Leitmann's method we prove it. For the classical (non-fractional)
approach, we refer the reader to \cite{silva,torres}.

\begin{Prop}
\label{prop:2}
Let $\alpha\in(0,1)$ be a real number. Consider
the following fractional optimal control problem:
\begin{equation}
\label{eq:ex2:funct}
\mathfrak{J}(u_1,u_2)
=\int_0^1\left[ (u_1(x))^2 +(u_2(x))^2 \right]\, dx \longrightarrow \min
\end{equation}
subject to the fractional control system
\begin{equation}
\label{system1}
\left\{
\begin{array}{l}
{_0D_x^\alpha}y_1(x)=\exp(u_1(x))+u_1(x)+u_2(x)\\
{_0D_x^\alpha}y_2(x)=u_2(x)\\
\end{array}\right.
\end{equation}
and fractional boundary conditions
\begin{equation}
\label{bound2}
\left\{
\begin{array}{ll}
{_0I_1^{1-\alpha}}y_1(1)=2\\
{_0I_1^{1-\alpha}}y_2(1)=1\, .\\
\end{array}\right.
\end{equation}
The global minimizer of problem
(\ref{eq:ex2:funct})--(\ref{bound2}) is
\begin{equation}
\label{sol:ex2}
(u_1(x),u_2(x))\equiv(0,1) \, , \quad
\left(y_1(x),y_2(x)\right)
=\left(\frac{2x^\alpha}{\alpha\Gamma(\alpha)},
\frac{x^\alpha}{\alpha\Gamma(\alpha)}\right) \, .
\end{equation}
\end{Prop}

\begin{proof}
We use the following change of variables:
$$
\left\{\begin{array}{l}
\tilde{y}_1(x)=y_1(x)-{_0I_x^{\alpha}}1\\
\tilde{y}_2(x)=y_2(x)-{_0I_x^{\alpha}}1\\
\tilde{u}_1(x)=u_1(x)\\
\tilde{u}_2(x)=u_2(x)-1\, .\\
\end{array}\right.
$$
For the new variables, we compute
$$\begin{array}{ll}
{_0D_x^\alpha}\tilde{y}_1& = {_0D_x^\alpha}y_1-{_0D_x^\alpha} {_0I_x^{\alpha}}1\\
&= \exp(u_1)+u_1+u_2-1\\
&= \exp(\tilde u_1)+\tilde u_1+\tilde u_2
\end{array}
$$
and
$$
\begin{array}{ll}
{_0D_x^\alpha}\tilde{y}_2
&= {_0D_x^\alpha}y_2-{_0D_x^\alpha} {_0I_x^{\alpha}}1\\
&=\tilde u_2.
\end{array}
$$
Therefore, $(\tilde y_1,\tilde y_2, \tilde u_1,\tilde u_2)$
satisfies (\ref{system1}). System (\ref{bound2}) becomes
\begin{equation}
\label{bound1}
{_0I_1^{1-\alpha}} \tilde y_1(1)=1
\quad \mbox{and} \quad {_0I_1^{1-\alpha}} \tilde y_2(1)=0
\end{equation}
since ${_0I_x^{1-\alpha}}{_0I_x^{\alpha}}1={_0I_x^{1}}1=x$.

Consider a new problem, which we label as $(\tilde P)$:
minimize $\mathfrak{J}(\cdot,\cdot)$ subject to system
(\ref{system1}) and conditions (\ref{bound1}).
Evaluating the functional $\mathfrak{J}(\cdot,\cdot)$
at $(\tilde u_1, \tilde u_2)$, we obtain
$$
\begin{array}{ll}
\mathfrak{J}(\tilde u_1,\tilde u_2)
&=\displaystyle\int_0^1\left[ (\tilde u_1(x))^2 +(\tilde u_2(x))^2 \right]\, dx\\
&=\displaystyle\int_0^1\left[ (u_1(x))^2 +(u_2(x))^2-2u_2(x)+1\right]\, dx \\
&=\mathfrak{J}(u_1,u_2)+\displaystyle\int_0^1 (-2{_0D_x^\alpha} y_2(x)+1)\, dx\\
&=\mathfrak{J}(u_1,u_2)+\displaystyle\int_0^1\frac{d}{dx} (-2 {_0I_x^{1-\alpha}} y_2(x) +x)\, dx.
\end{array}
$$
Observe that the second term of the last expression is constant.
Because $\tilde{u}_1(x)\equiv\tilde{u}_2(x)\equiv0$
is a solution of the problem $(\tilde P)$, we get
$$
\left\{
\begin{array}{l}
{_0D_x^\alpha} \tilde y_1=1\\
{_0D_x^\alpha} \tilde y_2=0\, ,\\
\end{array}\right.
$$
\textrm{i.e.},
$$
\left\{
\begin{array}{l}
\tilde y_1(x)={_0I_x^{\alpha}}1=\frac{x^\alpha}{\alpha\Gamma(\alpha)}\\
\tilde y_2(x)={_0I_x^{\alpha}}0=0 \, .
\end{array}\right.
$$
Therefore, the solution of problem
(\ref{eq:ex2:funct})--(\ref{bound2})
is given by
$$
(u_1(x),u_2(x))=(0,1)
$$
and
$$
\left(y_1(x),y_2(x)\right)
=\left(\frac{2x^\alpha}{\alpha\Gamma(\alpha)},
\frac{x^\alpha}{\alpha\Gamma(\alpha)}\right).
$$
\end{proof}

For $\alpha = 1$ the problem (\ref{eq:ex2:funct})--(\ref{bound2})
was solved in \cite[pag.~56]{torres}: the solution coincides
with (\ref{sol:ex2}) taking $\alpha=1$. It is worth noticing
that it is far from trivial to obtain Proposition~\ref{prop:2}
using the techniques of \cite{AGRA2,Frederico:Torres:2008}.


\section{Conclusions}
\label{sec:conc}

In this paper we illustrate the direct approach
of George Leitmann within the fractional variational context.
Although we only consider the continuous
fractional calculus of variations, the techniques
remain valid for the discrete variational calculus
recently introduced in \cite{comNunoRui:confOrtigueira09}.

To the best of our knowledge,
no direct methods have ever been applied before
to fractional optimal control or variational problems.
In our opinion there is an urgent need to promote
such direct methods among the community
of fractional variational problems.
Indeed, as illustrated in Section~\ref{sec4},
and notwithstanding the simplicity of
the statement of Theorem~\ref{thm:mr},
in some cases a difficult problem
can be simplified and easily solved analytically without
numerical solution of complicated fractional Euler-Lagrange type conditions.
More than that, Leitmann's approach allows one to solve problems
for which available methods of the literature
of fractional calculus fail to apply
(\textrm{cf.} Remark~\ref{rem:resp:r2}).


\section*{Acknowledgements}

Work partially supported by the \emph{control theory group}
(cotg) of the \emph{Centre for Research on Optimization
and Control} (CEOC) at the University of Aveiro,
with funds from \emph{The Portuguese Foundation for Science
and Technology} (FCT), cofinanced
by the European Community Fund FEDER/POCI 2010.

\smallskip

We are very grateful to three anonymous referees.



\end{document}